\documentclass{amsart}
\usepackage{amssymb}
\begin{document}
\newtheorem{theorem}{Theorem}[section]
\newtheorem{lemma}[theorem]{Lemma}
\newtheorem{corollary}[theorem]{Corollary}
\newtheorem{remark}[theorem]{Remark}
\def\operatorname#1{{\rm#1\,}}
\def\text#1{{\hbox{#1}}}
\def\qedbox{\hbox{$\rlap{$\sqcap$}\sqcup$}}
\def\range{\operatorname{range}}
\def\Pspan{\operatorname{span}}
\def\rank{\operatorname{rank}}
\def\id{\operatorname{Id}}
\def\trace{\operatorname{trace}}
\def\imag{\operatorname{Im}}
\font\pbglie=eufm10
\def\RR{\text{\pbglie R}}
\def\SS{\text{\pbglie S}}
\def\PP{\text{\pbglie P}}
\def\BB{A}
\def\Pspec{\operatorname{Spec}}
 \makeatletter
  \renewcommand{\theequation}{%
   \thesection.\arabic{equation}}
  \@addtoreset{equation}{section}
 \makeatother
\title[Szab\'o algebraic curvature tensors]{Jordan Szab\'o algebraic covariant derivative curvature tensors}
\author{Peter B. Gilkey, Raina Ivanova, and Iva Stavrov}
\thanks{Research of PG partially supported by the NSF (USA) and the MPI (Leipzig)}
\thanks{Research of RI and IS partially supported by the NSF (USA)}
\begin{address}{PG: Max Planck Institute in the Mathematical Sciences, 22 Inselstrasse, 04103 Leipzig Deutschland and
Mathematics Department, University of Oregon, Eugene OR 97403 USA}
\end{address}
\begin{email}
{gilkey@darkwing.uoregon.edu}
\end{email}
\begin{address}{IS: Mathematics Department, University of Oregon, Eugene OR 97403 USA}
\end{address}
\begin{email}{stavrov@hopf.uoregon.edu}\end{email}
\begin{address}{RI: Mathematics Department,
University of Hawaii - Hilo,
200 W. Kawili St.,
Hilo,  HI 96720 USA}\end{address}
\begin{email}
{rivanova@hawaii.edu}
\end{email}
\begin{abstract} We show that if
$\RR$ is a  Jordan Szab\'o algebraic covariant derivative curvature tensor on a vector space
of signature $(p,q)$, where $q\equiv1$ mod $2$ and $p<q$ or $q\equiv2$ mod $4$ and $p<q-1$, then $\RR=0$. This
algebraic result yields an elementary proof of the geometrical fact that any pointwise totally isotropic pseudo-Riemannian
manifold with such a signature
$(p,q)$ is locally symmetric.
\end{abstract}
\maketitle
\def\AA{{\mathcal{A}}}
\def\CC{{\mathcal{C}}}
\def\JJ{{\mathcal{J}}}
\def\OO{{\mathcal{O}}}
\font\pbglie=eufm10
\def\pp{\phantom{-}}
\def\Gr{\text{Gr}}
\section{Introduction}
Let $V$ be a vector space which is equipped with a non-degenerate
inner product $(\cdot,\cdot)$ of signature
$(p,q)$. Let $m:=p+q=\dim\{V\}$. Let
$S^\pm(V)$ be the pseudo-spheres of unit spacelike ($+$) and timelike ($-$) vectors in $V$:
$$S^\pm(V):=\{v\in V:(v,v)=\pm1\}.$$

We now review briefly some facts from linear algebra which we shall need in what follows. Let $A$ be a linear map from $V$ to
$V$. We say
$A$ is {\it self-adjoint} if
$(Av,w)=(v,Aw)$ for all
$v,w\in V$; let $\mathcal{A}(V)$ be the vector space of all self-adjoint linear maps from $V$ to $V$.
If $A$ is self-adjoint and if $\lambda\in\mathbb{C}$, then define real
operators
$A_\lambda$ on $V$ and associated generalized subspaces $E_\lambda$ by setting:
\begin{eqnarray}
  &&A_\lambda:=\left\{\begin{array}{ll}
A-\lambda\cdot\id&\quad\text{if}\quad\lambda\in\mathbb{R},\\
(A-\lambda\cdot\id)(A-\bar\lambda\cdot\id)&\quad\text{if}\quad\lambda\in\mathbb{C}-\mathbb{R},\end{array}\right.\nonumber\\
  &&E_\lambda:=\ker\{A_\lambda^m\}.\label{eqn1.1}\end{eqnarray}
Since $A_\lambda=A_{\bar\lambda}$, $E_\lambda=E_{\bar\lambda}$. Both $A$ and $A_\lambda$ preserve each
generalized eigenspace $E_\lambda$. The operator $A$
is said to be {\it Jordan simple} if $A_\lambda=0$ on $E_\lambda$ for all $\lambda$. We say $\lambda$ is an {\it eigenvalue
of $A$} if $\dim\{E_\lambda\}>0$; the {\it spectrum} $\Pspec\{A\}\subset\mathbb{C}$ is the set
of complex eigenvalues of
$A$.

\begin{lemma}\label{lem1.1}
Let $A$ be a self-adjoint linear map of a vector space of signature $(p,q)$.
\begin{enumerate}
\item We may decompose
$V=\oplus_{\imag(\lambda)\ge0}E_\lambda$.
\item We have $E_\lambda\perp E_\mu$
if $\lambda\ne\mu$ and $\lambda\ne\bar\mu$.
\item  The spaces $E_\lambda$ inherit non-degenerate metrics of signature
$(p_\lambda,q_\lambda)$.
\item If
$\lambda\notin\mathbb{R}$, then $p_\lambda=q_\lambda$.
\end{enumerate}\end{lemma} 

We say that $R\in\otimes^4V^*$ is an {\it algebraic curvature tensor} on $V$ if $R$ satisfies the symmetries:
\begin{eqnarray}
&&R(x,y,z,w)=R(z,w,x,y)=-R(y,x,z,w),\nonumber\\
&&R(x,y,z,w)+R(y,z,x,w)+R(z,x,y,w)=0.\label{eqn1.2}\end{eqnarray}
We say that
$\RR\in\otimes^5V^*$ is an {\it algebraic covariant derivative curvature tensor} on $V$ if $\RR$ satisfies the
symmetries:
\begin{eqnarray}
&&\RR(x,y,z,w;v)=\RR(z,w,x,y;v)=-\RR(y,x,z,w;v),\nonumber\\
&&\RR(x,y,z,w;v)+\RR(y,z,x,w;v)+\RR(z,x,y,w;v)=0,\label{eqn1.3}\\
&&\RR(x,y,z,w;v)+\RR(x,y,w,v;z)+\RR(x,y,v,z;w)=0.\nonumber
\end{eqnarray}

Let $(M,g)$ be a pseudo-Riemannian manifold of signature $(p,q)$. The Riemann curvature tensor $R_g$ is
an algebraic curvature tensor on the tangent space
$T_PM$ for every point $P\in M$. Similarly, the covariant derivative of the curvature tensor $\nabla R_g$ is an
algebraic covariant derivative curvature tensor on $T_PM$. Conversely, such tensors are geometrically realizable. Thus
tensors with the symmetries of equations (\ref{eqn1.2}) and (\ref{eqn1.3}) are important in differential geometry.

Let $R$ be an algebraic curvature tensor on $V$. The associated {\it Jacobi operator} $J$ is the linear map of
$V$ characterized by the identity:
$$(J(v)y,z)=R(y,v,v,z);$$
$J(v)$ is self-adjoint by equation (\ref{eqn1.2}). We say that $R$ is a {\it spacelike
Osserman} (resp. {\it timelike Osserman}) tensor if
$R$ is an algebraic curvature tensor and if $\Pspec\{J\}$ is constant on $S^+(V)$ (resp. on $S^-(V)$); these two notions are
equivalent
\cite{refGil} and such a tensor is said to be an {\it Osserman} tensor.

Analogously, let $\RR$ be an algebraic covariant derivative curvature tensor on $V$.
The associated {\it Szab\'o operator} $\SS$ is the linear map of $V$ characterized by
the identity:
$$(\SS(v)y,z)=\RR(y,v,v,z;v);$$
the symmetries of equation (\ref{eqn1.3}) show $\SS$ is self-adjoint. We say
that
$\RR$ is a {\it spacelike Szab\'o} (resp. {\it timelike Szab\'o}) tensor if $\RR$ is an algebraic covariant derivative
curvature tensor and if
$\Pspec\{\SS\}$ is constant on
$S^+(V)$ (resp. $S^-(V)$); we denote this common spectrum by $\Pspec^+\{\SS\}$ (resp. $\Pspec^-\{\SS\}$). 
The notions spacelike Szab\'o and timelike Szab\'o are equivalent \cite{refGil}; such a tensor is said to be a {\it
Szab\'o tensor}.

Osserman \cite{refOss} and Szab\'o \cite{refSzabo} wrote the original papers concerning the
spectral properties of the operators $J$ and $\SS$ in the Riemannian
setting. That is why their names have become associated with the subject. The Jacobi operator has been studied
extensively in this context; we refer to
\cite{refGRKVL,refGil} for a more complete bibliography. By contrast, the Szab\'o operator has received considerably
less attention and the present paper is devoted to the study of the spectral properties of $\SS$ in the
pseudo-Riemannian setting. Although there are certain formal parallels between the Jacobi operator and the Szab\'o
operator, the fact that
$J(-v)=J(v)$ while
$\SS(v)=-\SS(-v)$ plays a crucial role.

In the higher signature setting, the eigenvalue structure does not control the Jordan normal form (i.e. the conjugacy
class) of a self-adoint map. We shall say that $\RR$ is a {\it spacelike} (resp. {\it timelike}) {\it Jordan Szab\'o
tensor} if
$\RR$ is an algebraic covariant derivative curvature tensor and if  the Jordan normal form of $\SS$ is
constant on $S^+(V)$ (resp. $S^-(V)$).

The following result is due in the Riemannian setting ($p=0)$ to Szab\'o \cite{refSzabo} and in the Lorentzian setting
($p=1$) to Gilkey and Stavrov \cite{refGilStav}:

\begin{theorem}\label{thm1.2} Let $\RR$ be a Szab\'o tensor on a vector space of signature $(p,q)$. If $p=0$ or if $p=1$,
then
$\RR=0$.\end{theorem}

By replacing $g$ by $-g$, one can interchange the roles of $p$ and $q$. Consequently Theorem \ref{thm1.2}, implies there
are no non-trivial Szab\'o tensors if $q=0$ or if $q=1$ as well. Therefore, we shall assume
$p\ge2$ and
$q\ge2$ henceforth. Although we are primarily interested in spacelike or timelike Jordan Szab\'o tensors,
many of our results hold under the weaker assumption that the tensor is Szab\'o. 
We shall
establish the following result in Section
\ref{Sect2}.

\begin{theorem}\label{thm1.3}  Let $\RR$ be
a Szab\'o tensor on a vector space $V$ of signature $(p,q)$, where $p\ge2$ and $q\ge2$.
Then\begin{enumerate}
\item $\Pspec^\pm\{\SS\}=-\Pspec^\pm\{\SS\}=\sqrt{-1}\phantom{.}\Pspec^\mp\{\SS\}$;
\item $\Pspec^\pm\{\SS\}\subset\mathbb{R}\cup\sqrt{-1}\phantom{.}\mathbb{R}$;
\item If $p<q$, then $\Pspec^+\{\SS\}\subset\sqrt{-1}\phantom{.}\mathbb{R}$ and $\Pspec^-\{\SS\}\subset\mathbb{R}$;
\item If $q<p$, then $\Pspec^+\{\SS\}\subset\mathbb{R}$ and $\Pspec^-\{\SS\}\subset\sqrt{-1}\phantom{.}\mathbb{R}$.
\end{enumerate}
\end{theorem}

A pseudo-Riemannian manifold $(M,g)$ of signature $(p,q)$ is said to be a {\it Szab\'o} manifold if
$\nabla R_g$ is a Szab\'o tensor on $T_PM$ for all $P\in M$. If $p\ge2$ and if $q\ge2$, then
 there exist Szab\'o pseudo-Riemannian manifolds and algebraic covariant derivative curvature tensors of signature $(p,q)$
which are neither locally symmetric, nor locally homogeneous, nor pointwise totally isotropic \cite{refFG02,refGIZ}.

Let $S^{q-1}$ be the unit sphere in $\mathbb{R}^q$ with the usual  positive definite inner product.
If $q$ is odd, then there are no non-vanishing vector fields on $S^{q-1}$. This well known result has been generalized
by Adams \cite{refAdams}. The {\it Adams number} $\nu(q)$ is defined as follows. Let $q=2^\ell q_0$, where $q_0$ is odd. Then
$\nu(q)=\nu(2^\ell)$, where
\begin{equation}\nu(1)=0,\quad \nu(2)=1,\quad \nu(4)=3,\ 
  \nu(8)=7,\quad\text{and}\quad\nu(2^{\ell+4})=\nu(2^\ell)+8.\label{eqn1.4}\end{equation}
Let $\{v_1,...,v_\ell\}$ be linearly independent vector fields on $S^{q-1}$. Adams showed that
$\ell\le\nu(q)$; thus $\nu(q)$ provides an upper bound to the number of linearly independent vector fields that can exist on
$S^{q-1}$. Furthermore, the estimate is sharp. If $\nu(q)>0$, then there exist $\nu(q)$ linearly independent vector fields
on $S^{q-1}$. Since $\nu(q)=0$ if $q$ is odd, the result of Adams contains as a special case the original observation that
there does not exist a nowhere vanishing vector field on an even dimensional sphere.

If $p\ge2$ and if $q\ge2$ , then there non-trivial nilpotent Szab\'o pseudo-Riemannian manifolds \cite{refFG02,refGIZ}.
However, there are as yet no known examples of spacelike Jordan Szab\'o tensors $\RR$, where
$\RR\ne0$, and we conjecture none exist. We have some partial results in this direction. In section \ref{Sect3}, we shall establish the following result:

\begin{theorem}\label{thm1.4} Let $\RR$ be a spacelike Jordan Szab\'o tensor on a vector space $V$ of signature
$(p,q)$, where $p<q$. Let $v\in S^+(V)$. Then:
\begin{enumerate}
\item $\SS(v)$ is Jordan simple;
\item If $p<q-\nu(q)$, then $\rank\{\SS(v)\}\le2\nu(q)$;
\item If $q$ is odd, then $\RR=0$.
\end{enumerate}
\end{theorem}

If $\RR$ is spacelike and timelike
Jordan Szab\'o, then we let $r_\pm$ be the rank of $\SS$ on the pseudo-spheres $S^\pm(V)$. In Section \ref{Sect4}, we shall
prove:

\begin{theorem}\label{thm1.5}  Let $\RR\ne0$ be an algebraic covariant derivative curvature tensor on a vector space $V$ of
signature
$(p,q)$ which is both spacelike and timelike Jordan Szab\'o. Then:
\begin{enumerate}
\item $r_+=r_-$.
\item If $p\ne q$, then $\SS(v)$ is Jordan simple for $v\in S^\pm(V)$.
\end{enumerate}
\end{theorem}

We say $\RR$ is {\it null Jordan Szab\'o} if the Jordan normal form of $\SS$ is constant on the null cone
$$\mathcal{N}:=\big\{\{v\in V:(v,v)=0\}-\{0\}\big\};$$
this implies $\SS$ is nilpotent and has constant rank, which we shall denote by $r_0$, on $\mathcal{N}$. We say that $\RR$ is
{\it Jordan Szab\'o} if $\RR$ is spacelike, timelike, and null Jordan Szab\'o. In Section \ref{Sect5}, we shall prove:
\begin{theorem}\label{thm1.6}  Let $\RR$ be an algebraic covariant derivative curvature tensor on a vector space $V$ of
signature
$(p,q)$ which is Jordan Szab\'o. Then
\begin{enumerate}
\item If $\RR\ne0$, then $r_0<r_+$.
\item If $q\equiv2$ mod $4$ and if $p<q-1$, then $\RR=0$.
\end{enumerate}\end{theorem}

Following Wolf \cite{refWolf}, we say that a pseudo-Riemannian manifold $(M,g)$ is {\it locally
isotropic} if given a point $P\in M$ and nonzero tangent vectors $X$ and $Y$ in $T_PM$ with
$(X,X)=(Y,Y)$, there is a local isometry of
$(M,g)$, fixing $P$, which sends $X$ to $Y$. Wolf
showed that such a manifold is necessarily locally symmetric, see \cite{refWolf} (Theorem 12.3.1).  The Szab\'o
operator of a locally isotropic pseudo-Riemannian manifold is necessarily Jordan Szab\'o. Thus
Wolf's result in certain special cases can be derived from Theorems \ref{thm1.4} and
\ref{thm1.5}:

\begin{corollary}\label{cor1.7} Let $(M,g)$ be a locally isotropic pseudo-Riemannian manifold of
signature $(p,q)$. If $q\equiv1$ mod $2$ and if $p<q$ or if $q\equiv2$ mod $4$ and if $p<q-1$, then $(M,g)$ locally
symmetric.
\end{corollary}

\section{Spacelike Szab\'o Tensors}\label{Sect2}

In this section we prove Theorem \ref{thm1.3}. Let $V$ be a vector space which is equipped with a non-degenerate
inner product of signature $(p,q)$, where $p\ge2$ and $q\ge2$. We can choose a non-canonical decomposition
\begin{equation}V=V^+\oplus V^-,\label{eqn2.1}\end{equation}
 where $V^+$ is a maximal
spacelike subspace and where $V^-:=(V^+)^\perp$ is the complementary maximal timelike subspace. Let $\pi^\pm$ be orthogonal
projection on
$V^\pm$. The following is a useful technical observation.
\begin{lemma}\label{lem2.1} Let $v\in S^+(V^+)$. If $w\perp v$, then $\pi^+w\in V^+\cap v^\perp=T_v(S^+(V^+))$.
\end{lemma}

\begin{proof} We may decompose $w=cv+v^++v^-$, where $v^\pm\in V^\pm$ and $v\perp v^+$. Because $c=(v,w)=0$, we have
$\pi^+w=v^+$.\end{proof}

The next result is well known, see, for example 
\cite{refGiIv}, and generalizes the decomposition of equation (\ref{eqn2.1}) to the bundle setting:

\begin{lemma}\label{lem2.2} Let $E$ be a vector bundle over a smooth manifold $M$ which is
equipped with a non-degenerate fiber metric. Then we can decompose $E$ as a direct sum
$E^+\oplus E^-$ of complementary orthogonal subbundles, where $E^+$ is a maximal spacelike subbundle and
$E^-=(E^+)^\perp$ is a maximal timelike subbundle.\end{lemma}

 Let $\RR$ be a Szab\'o tensor on $V$.
Since $\RR(-v)=-\RR(v)$,
$$\Pspec\{\SS(v)\}=-\Pspec\{\SS(-v)\}\quad\text{so}\quad\Pspec^\pm\{\SS\}=-\Pspec^\pm\{\SS\}.$$

We now establish Theorem \ref{thm1.3} (1). Let $\RR$ be a Szab\'o tensor on a vector space of signature $(p,q)$ for $p\ge2$
and $q\ge2$. Let
$V_{\mathbb{C}}:=V\otimes\mathbb{C}$ be the complexification of $V$. We extend the inner product and the tensor $\RR$ to be
complex multilinear on $V_{\mathbb{C}}$. Let
\begin{eqnarray*}
\mathcal{O}_{\mathbb{C}}:&=&\{v\in V_{\mathbb{C}}:(v,v)\ne0\},\\
\PP(v):&=&(v,v)^{-3}\SS(v)^2\text{ for }v\in\mathcal{O}_{\mathbb{C}},\text{ and}\\
p(t,v):&=&\det\{\PP(v)-t\cdot\id\}=c_0(v)+...+c_m(v)t^m\text{ for }v\in\mathcal{O}_{\mathbb{C}}.\end{eqnarray*}
The set $\mathcal{O}_{\mathbb{C}}$ is a connected open subset of $V_{\mathbb{C}}$ and the coefficients $c_i(v)$ are
holomorphic functions on $\mathcal{O}_{\mathbb{C}}$. The intersection of $\mathcal{O}_{\mathbb{C}}$ with the real subspace
$V$ is the disjoint union of the non-empty sets $\mathcal{O}_\pm$ of spacelike and timelike vectors in $V$. Since
$\RR$ is a Szab\'o tensor, $\Pspec\{\PP\}$, the characteristic polynomial $p$, and hence the
functions
$c_i$ are constant on $\mathcal{O}_+$. By the identity theorem, the functions $c_i$ are constant on
$\mathcal{O}_{\mathbb{C}}$. Let $v_\pm\in S^\pm(V)$. Because
$p(t,v^-)=p(t,v^+)$, the operators
$\PP(v^-)=-\SS(v^-)^2$ and $\PP(v^+)=\SS(v^+)^2$ have the same spectrum. This completes the proof of assertion (1) of
Theorem
\ref{thm1.3} by establishing the identity:
$$\Pspec\{\SS(v^-)\}=\sqrt{-1}\phantom{.}\Pspec\{\SS(v^+)\}.$$

Before continuing to prove the other assertions of Theorem \ref{thm1.3}, we must recall some facts from algebraic topology.
Assertion (1) in the following Lemma follows from results of Szab\'o
\cite{refSzabo}, and assertion (2) follows from the Borsuk-Ulam theorem (see for example Spanier
\cite{refSpanier}, Corollary 8 p. 266). They concern $\mathbb{Z}_2$ equivariant vector fields and vector valued functions
on spheres.

\begin{lemma}\label{lem2.3} Let
$S^{q-1}$ be the unit sphere in a vector space of signature $(0,q)$. 
\begin{enumerate}
\smallskip\item There does not exist a continuous nowhere vanishing vector field $F$ on $S^{q-1}$ so $F(v)=F(-v)$ for every
$v\in V$.
\smallskip\item If $\dim\{V^-\}=p<q$, then there does not exist a continuous nowhere vanishing map $F$ from
$S^{q-1}$ to $V^-$ so that
$F(v)=-F(-v)$ for every
$v\in S^{q-1}$.
\end{enumerate}
\end{lemma}

Recall that $\mathcal{A}(V)$ is the set of all self-adjoint linear maps of $V$. We adopt the notation of Lemma \ref{lem1.1}.
Let $V=V^+\oplus V^-$ be the decomposition described in equation (\ref{eqn2.1}).
The remaining assertions of Theorem \ref{thm1.3} will be proved, as has become traditional in this subject, using
methods from algebraic topology. The essential point is to abstract precisely those properties
relevant to the investigation. What will be
crucial for us are the facts that:
\begin{eqnarray}
&&\SS(-v)=-\SS(v),\quad
\SS(v)v=0,\quad\text{and}\nonumber\\
&&\SS(v)\in\mathcal{A}(V)\quad\text{for all}\quad v\in
S^{q-1}:=S^+(V^+).\label{eqn2.2}\end{eqnarray}

\begin{lemma}\label{lem2.4} Let $V$ be a vector space of signature $(p,q)$. Let $V^+$ be a maximal spacelike
subspace of $V$ and let $S^{q-1}:=S^+(V^+)$. Let $A:v\rightarrow A_v$ be a continuous map from
$S^{q-1}$ to $\mathcal{A}(V)$ with $A_{-v}=-A_v$ and with $A_vv=0$ for all $v\in S^{q-1}$.
Assume that $A$ has constant spectrum. Then:
\begin{enumerate}\item If $\lambda\in\Pspec\{A\}$ with $\lambda\notin\sqrt{-1}\phantom{.}\mathbb{R}$, then $q_\lambda=0$;
\item $\Pspec\{A\}\subset\mathbb{R}\cup\sqrt{-1}\phantom{.}\mathbb{R}$;
\item If $p<q$, then $\Pspec\{A\}\subset\sqrt{-1}\phantom{.}\mathbb{R}$.\end{enumerate}\end{lemma}

\begin{proof} Let $\pi^\pm$ be orthogonal
projection on $V^\pm$. Consider the following vector bundles on
$S^{q-1}$:
\begin{eqnarray}
TS^{q-1}&:=&\{(v_1,v_2)\in S^{q-1}\times V^+:(v_1,v_2)=0\},\text{ and}\nonumber\\
TS^{q-1}\oplus V^-&:=&\{(v_1,v_2)\in S^{q-1}\times V:(v_1,v_2)=0\}.\nonumber\end{eqnarray}
Let $A:S^{q-1}\rightarrow\mathcal{A}(V)$ satisfy the hypothesis of the Lemma. Since $A$ has constant spectrum, the
generalized eigenspaces $E_{\lambda,v}$ of $A_v$, as described in Lemma \ref{lem1.1}, have constant rank and fit together to
define smooth vector bundles
$E_\lambda$ over
$S^{q-1}$.

Suppose $\lambda\in\Pspec\{A\}$ and $\lambda\notin\sqrt{-1}\phantom{.}\mathbb{R}$. We apply Lemma
\ref{lem1.1}; $0\ne\lambda$, $-\lambda\ne\lambda$ and $-\lambda\ne\bar\lambda$. Thus $E_{0,v}$, $E_{\lambda,v}$, and
$E_{-\lambda,v}$ are mutually orthogonal vector spaces which inherit non-degenerate metrics of signatures $(p_0,q_0)$,
$(p_\lambda,q_\lambda)$, and $(p_{-\lambda},q_{-\lambda})$, respectively. As
$A_vv=0$,
$v\in E_{0,v}$ so $v\perp E_{\lambda,v}$. By Lemma \ref{lem2.2}, there exist
maximal orthogonal timelike and spacelike sub-bundles
$E_\lambda^\pm$ of $E_\lambda$ so:
$$E_\lambda=E^+_\lambda\oplus E^-_\lambda.$$
Since $A_{-v}=-A_v$, we have $E_{\lambda,-v}=E_{-\lambda,v}$. Thus we may obtain a similar
splitting of the bundle $E_{-\lambda}$ by setting:
\begin{eqnarray*}
&&E^\pm_{-\lambda,v}:=E^\pm_{\lambda,-v}.\end{eqnarray*}
Let $\varepsilon$ be a choice of signs $\pm$. We suppose $\dim\{E^\varepsilon_\lambda\}>0$; this means
$p_\lambda>0$ if $\varepsilon$ is $-$ and $q_\lambda>0$ if $\varepsilon$ is $+$. 
Let $N$ be the north pole of $S^{q-1}$. Since $S^{q-1}-\{N\}$ is contractable, there exists a global section $s_\lambda$ to
$E^\varepsilon_\lambda$ which only vanishes at $N$. Let $s_{-\lambda}(v):=s_{\lambda}(-v)$. Then $s_{-\lambda}$ is a
section to
$E^\varepsilon_{-\lambda}$ which only vanishes at the south pole $S:=-N$ and which satisfies $s_\lambda(v)\perp
s_{-\lambda}(v)$. Let
$$F:=s_\lambda+\varepsilon s_{-\lambda}.$$
If $\varepsilon$ is $+$, then $s_\lambda$ and $s_{-\lambda}$ are spacelike; if $\varepsilon$ is $-$, then $s_\lambda$ and
$s_{-\lambda}$ are timelike. Since $s_\lambda\perp s_{-\lambda}$, $F$ is a nowhere vanishing section to
$E^\varepsilon_\lambda\oplus E^\varepsilon_{-\lambda}$.

To prove assertion (1) of Lemma \ref{lem2.4}, we suppose that $q_\lambda>0$ and argue for a contradiction. We take
$\varepsilon$ to be
$+$. Then
$F(v)$ is nowhere vanishing, spacelike, and perpendicular to $v$ so Lemma \ref{lem2.1} shows that
$\pi^+F(v)$ is a nowhere vanishing vector field on $S^{q-1}$ with the equivariance property $\pi^+F(v)=\pi^+F(-v)$. This
contradicts Lemma
\ref{lem2.3} (1) and proves Lemma \ref{lem2.4} (1). 

\medbreak If $\lambda\in\Pspec\{A\}$ and if $\lambda\notin\mathbb{R}$, then $p_\lambda=q_\lambda\ne0$. Assertion (1)
then implies $\lambda\in\sqrt{-1}\phantom{.}\mathbb{R}$. Thus $\Pspec\{A\}\subset\mathbb{R}\cup\sqrt{-1}\phantom{.}\mathbb{R}$
which establishes assertion (2).

\medbreak Let $p<q$. To prove assertion (3), we assume that there exists an eigenvalue $\lambda$ belonging to
$\mathbb{R}-\{0\}$ and argue for a contradiction. We take
$\varepsilon$ equal to $-$. As
$F(v)$ is timelike and non-vanishing,
$\pi^-F(v)$ is a nowhere vanishing continuous map from $S^{q-1}$ to $V^-$ with the equivariance
property $\pi^-F(v)=-\pi^-F(-v)$. This contradicts Lemma
\ref{lem2.3} (2) and proves assertion (3) of Lemma \ref{lem2.4}. 
\end{proof}

\begin{proof}[Proof of Theorem \ref{thm1.3}] Let $\RR$ be a
spacelike Szab\'o tensor on
$V$. We have already established assertion (1) of the Theorem. Assertions (2), (3), and (4) follow from Lemma \ref{lem2.4},
where, if necessary, we change the sign of the quadratic form to interchange the roles of $p$ and $q$ and the notions
spacelike and timelike.
\end{proof}

\section{Spacelike Jordan Szab\'o Tensors}\label{Sect3}

We begin the preparation of Theorem \ref{thm1.4} by recalling some topological results. Let
$\mathbb{RP}^{q-1}$ be real projective space; this is the space of lines thru the origin in $\mathbb{R}^q$. If $v\in
S^{q-1}$, then we denote the corresponding element of $\mathbb{RP}^{q-1}$ by $\langle v\rangle=v\cdot\mathbb{R}$. We
may identify
$\mathbb{RP}^{q-1}=S^{q-1}/\mathbb{Z}_2$. The  classifying line bundle and orthogonal complement bundle over projective space
are defined by
\begin{eqnarray*}
L&:=&\{(\langle v\rangle,w)\in\mathbb{RP}^{q-1}\times\mathbb{R}^q:w\in\langle v\rangle\}\\
L^\perp&:=&\{(\langle v\rangle,w)\in\mathbb{RP}^{q-1}\times\mathbb{R}^q:v\perp w\}.\end{eqnarray*}
Note that we may identify $T(\mathbb{RP}^{q-1})$ with $L\otimes L^\perp$.
Let $\nu(q)$ be the Adams number which was defined in equation (\ref{eqn1.4}). We refer to \cite{refGiIv} for the proof of
the first assertion and to Adams \cite{refAdams} for the proof of the second assertion in the following Lemma:

\begin{lemma}\label{lem3.1}\ \begin{enumerate}
\item Let $U_i$ be non-trivial vector bundles over $\mathbb{RP}^{q-1}$. If $U_1$ is
a sub-bundle of
$L^\perp$ and if $U_2$ is a sub-bundle of a trivial bundle of dimension $p<q$, then $U_1$ is not isomorphic to $U_2$.
\item Let $S^{q-1}$ be the unit sphere in $\mathbb{R}^q$. Let $\nu(q)$ be the Adams number. Then:
\begin{enumerate}\item If $E$ is a rank $r$ sub-bundle of $TS^{q-1}$, then either
$q-1-\nu(q)\le r$ or $r\le\nu(q)$;
\item If $\{\chi_1,...,\chi_\ell\}$ are linearly independent vector fields on $S^{q-1}$, then we have $\ell\le\nu(q)$.
\end{enumerate}
\end{enumerate}\end{lemma}

Again, the properties of equation (\ref{eqn2.2}) are crucial. We begin the proof of Theorem
\ref{thm1.4} by using Lemma \ref{lem3.1} to establish a related result in the abstract setting. As in Lemma \ref{lem2.4}, we
decompose $V=V^-\oplus V^+$ and identify $S^{q-1}$ with the unit sphere in $V^+$.

\begin{lemma}\label{lem3.2}  Let $V$ be a vector space of signature $(p,q)$, where $p<q$. Let $V^+$ be a maximal spacelike
subspace of $V$ and let $S^{q-1}:=S^+(V^+)$. Let $A:v\rightarrow A_v$ be a continuous map from
$S^{q-1}$ to $\mathcal{A}(V)$ with $A_{-v}=-A_v$ and with $A_vv=0$ for all $v\in S^{q-1}$.
Assume that $A$ has constant Jordan normal form. Then:
\begin{enumerate}
\item The operator $A_v$ is Jordan simple for any $v\in S^{q-1}$;
\item If $p<q-\nu(q)$, then $\rank\{A)\}\le2\cdot\nu(q)$.
\item If $q$ is odd, then $A=0$.
\end{enumerate}
\end{lemma}

\begin{proof} By Lemma \ref{lem2.4},
$\Pspec\{A\}\subset\sqrt{-1}\phantom{.}\mathbb{R}$. We use equation (\ref{eqn1.1}) to define the operators $A_{\lambda,v}$
for $\lambda\in\Pspec\{A_v\}$. Since $\lambda$ is purely imaginary,
\begin{equation}A_{\lambda,v}=\left\{\begin{array}{ll}
A_v^2+|\lambda|^2&\text{if}\quad\lambda\ne0,\\
A_v&\text{if}\quad\lambda=0.\end{array}\right.\label{eqn3.1}\end{equation}
Since $A_vv=0$, $A_v$ preserves $v^\perp$. We therefore introduce the {\it reduced generalized eigenspaces}
$$\tilde E_{v,\lambda}:=\ker\{A_{\lambda,v}^m\}\cap v^\perp.$$

To prove assertion (1), we suppose, to the contrary, that $A_{\lambda,v}$ is not Jordan simple. Since $A_vv=0$, this means
that
$A_{\lambda,v}\ne0$ on the generalized eigenspaces 
$\tilde E_{v,\lambda}$.
Choose $k\ge1$ maximal so
$A_{\lambda,v}^k\ne0$ on $\tilde E_{\lambda,v}$. Then $A_{\lambda,v}^{k+1}=0$. Let 
$$U_{\lambda,v}:=A_{\lambda,v}^k\tilde E_{v,\lambda}.$$
Because $A_v$ has
constant Jordan normal form, the vector spaces $U_{\lambda,v}$ glue together smoothly to define a smooth vector bundle
$U_\lambda$ over
$S^{q-1}$. This bundle descends to a smooth bundle
$W_\lambda$ over $\mathbb{RP}^{q-1}=S^{q-1}/\mathbb{Z}_2$ because, by equation (\ref{eqn3.1})
$$A_{\lambda,v}=\pm A_{\lambda,-v}.$$

Let $w_i=A_{\lambda,v}^kv_i$ for $v_i\in\tilde E_{v,\lambda}$. Since $k\ge1$, $2k>k$ so:
\begin{equation}
0=(A_{\lambda,v}^{2k}v_1,v_2)=(A_{\lambda,v}^kv_1,A_{\lambda,v}^kv_2)=(w_1,w_2).\label{eqn3.2}\end{equation} 
By equation (\ref{eqn3.2}), the bundle $W_\lambda$ is totally isotropic.
Since $W_\lambda$ contains no spacelike or timelike vectors, the projections $\pi^\pm$ are non-singular
on
$W_\lambda$. Set 
$$W_\lambda^\pm:=\pi^\pm(W_\lambda).$$
If $w\in W_{\lambda,v}$, then $w\perp v$. Consequently, by Lemma \ref{lem2.1}, $\pi^+w\in V^+\cap v^\perp$. It is immediate
that
$\pi^-w\in V^-$. Therefore,
$$W_\lambda^+\subset L^\perp\quad\text{and}\quad W_\lambda^-\subset V^-.$$
As  $W_\lambda^+$, $W_\lambda$, and $W_\lambda^-$ are isomorphic, this
contradicts Lemma \ref{lem3.1} and proves Lemma \ref{lem3.2} (1).

By assertion (1), $A_v=0$ on $E_{0,v}$. Thus
\begin{equation}\rank\{A\}=2\sigma\quad\text{where}\quad
\sigma:=\textstyle\sum_{\Im(\lambda)>0}p_\lambda=\textstyle\sum_{\Im(\lambda)>0}q_\lambda.\label{eqn3.3}\end{equation}
By Lemma \ref{lem2.2}, choose a maximal spacelike sub-bundle $\tilde E_{v,\lambda}^+$  of
$\tilde E_{v,\lambda}$ of rank
$q_\lambda$. Let
\begin{equation}\tilde E^+=\oplus_{\Im(\lambda)>0}\tilde E_{v,\lambda}^+\label{eqn3.4}\end{equation}
be a sub-bundle of $v^\perp$ of rank $\sigma$. Since elements of $\tilde E^+$ are all spacelike, the projection $\pi^+$ is
injective and, by Lemma \ref{lem2.1}, $\tilde E^+$ is a sub-bundle of $TS^{q-1}$ of rank $\sigma$. By Lemma
\ref{lem3.1} (2), either
$\sigma\le\nu(q)$ or $\sigma\ge q-1-\nu(q)$. Since the first inequality implies assertion (2) of Lemma \ref{lem3.2}, we
assume, to the contrary, that
\begin{equation}\sigma\ge q-1-\nu(q).\label{eqn3.5}\end{equation}

Displays (\ref{eqn3.3}) and (\ref{eqn3.5}) imply 
\begin{equation}p\ge\sigma\ge q-1-\nu(q).\label{eqn3.6}\end{equation}
We assumed that $p<q-\nu(q)$. Therefore, by (\ref{eqn3.6}), we have:
$$\sigma=p=q-1-\nu(q).$$
Let $e\in\tilde E_{\lambda,v}^+$. By assertion (1), $A_v$ is Jordan simple. Thus $( A_v^2+|\lambda|^2)e=0$
so
\begin{equation}(A_ve,A_ve)=(A_v^2e,e)=-|\lambda|^2(e,e).\label{eqn3.7}\end{equation}
The decomposition of equation (\ref{eqn3.4}) is an orthogonal direct sum. Equation (\ref{eqn3.7}) shows that 
$A_v\tilde E^+_{\lambda,v}$ is a timelike subspace. Because $\rank\{A_v\tilde E^+_v\}=p$, $\pi^-$ is an isomorphism from
$A_v\tilde E^+$ onto
$V^-$. Consequently, $\pi^+\tilde E^+$ is isomorphic to a trivial bundle so there exist $p$ linearly independent vector
fields on
$S^{q-1}$. This shows $\sigma\le p\le\nu(q)$ and completes the proof of assertion (2).

Assertion (3) follows directly from assertion (2) since $\nu(q)=0$ if $q$ is odd.\end{proof}


\medbreak The following observation is well known, see, for example, \cite{refGilStav}:

\begin{lemma}\label{lem3.3} Let $\RR$ be an algebraic covariant derivative curvature tensor. If $\SS$ vanishes identically,
then $\RR=0$.\end{lemma}

\begin{proof}[Proof of Theorem \ref{thm1.4}] Let $\RR$ be a spacelike Jordan Szab\'o tensor on a vector space $V$ of
signature $(p,q)$, where $p<q$. By equation (\ref{eqn2.2}), we may apply Lemma \ref{lem3.2}. Assertions (1) and (2) of
Theorem \ref{thm1.4} now follow directly from assertions (1) and (2) of Lemma \ref{lem3.2}. If $q$ is odd, we apply Lemma
\ref{lem3.2} (3) to see
$\SS=0$. Assertion (3) now follows from Lemma \ref{lem3.3}. \end{proof}

\section{Timelike and spacelike Jordan Szab\'o Tensors}\label{Sect4}

Let $\RR$ be both timelike and spacelike Jordan Szab\'o. Let $r_\pm:=\rank\{\SS(\cdot)\}$ on $S^\pm(V)$. To establish Theorem
\ref{thm1.5} (1), we must show that
$r_+=r_-$. Fix a vector $v^+\in S^+(V)$. Choose vectors $\{v_1,...,v_{r_+}\}$ in $V$ and choose dual
elements
$\{v_1^*,...,v_{r_+}^*\}$ in the dual vector space $V^*$ so  
$$v_i^*\SS(v^+)v_j=\delta_{ij}.$$
We adopt the notation
established in the proof of Lemma
\ref{lem2.2} and complexify. If $v\in V_{\mathbb{C}}$, then we consider the matrix $A_{ij}(v):=v_i^*\SS(v)v_j$. We
define a holomorphic function on $V_{\mathbb{C}}$ by setting:
$$f(v):=\det\{v_i^*\SS(v)v_j\}.$$
Since $f(v^+)=1$, $f$ does not vanish identically on $V_{\mathbb{C}}$.
Thus by the identity theorem, there must exist a real timelike element $v^-$ of $V$ so $f(v^-)\ne0$. This implies that the
elements
$\{\SS(v^-)v_1,...,\SS(v^-)v_{r_+}\}$ are linearly independent so
$$r_-=\rank\{\SS(v^-)\}\ge r_+.$$
Similarly we can show $r_+\ge r_-$. Thus
$r_+=r_-$.

 To prove the second assertion, we may suppose without loss of generality that $p<q$. If $v\in V_{\mathbb{C}}$, then
we shall let
$E^c_\lambda\subset V_\mathbb{C}$ be the generalized complex eigenspaces of
$\SS(v)$. Let $v_\pm\in S^\pm(V)$. As, by Theorem \ref{thm1.4}, $\SS(v^+)$ is Jordan simple,
$$r_+=\textstyle\sum_{\lambda\ne0}\dim_{\mathbb{C}}\{E^c_\lambda\}.$$
In Section \ref{Sect2}, we showed that the characteristic polynomials of the two operators $\SS(v^+)$ and
$\sqrt{-1}\phantom{.}\SS(v^-)$ were the same. Thus
$$r_+=\textstyle\sum_{\lambda\ne0}\dim_{\mathbb{C}}\{E^c_\lambda(\SS(v^-))\}.$$
Since $r_+=r_-$, we conclude that $\SS(v^-)$ must have rank $0$ on $E^c_0(\SS(v^-))$ and consequently $0$ is a Jordan simple
eigenvalue for $\SS(v^-)$. By Theorem \ref{thm1.3}, the eigenvalues of $\SS(v^-)$ are real. By Lemma \ref{lem2.4}, the
generalized eigenspaces for $\lambda\in\mathbb{R}-\{0\}$ are spacelike. Since $\SS(v^-)$ is self-adjoint with respect to a
definite metric on $E_\lambda$, $\SS(v^-)$ is diagonalizable on $E_\lambda$. This shows $\SS(v^-)$ is Jordan simple and
establishes assertion (2) and completes the proof of Theorem \ref{thm1.5}.

\section{Jordan Szab\'o Tensors}\label{Sect5}

\begin{proof}[Proof of Theorem \ref{thm1.6} (1)] Let $\RR$ be a non-trivial Jordan Szab\'o tensor whose rank on the null
cone is $r_0$. We must show that
$r_0<r_+$. Let $v_0\in\mathcal{N}$. As in the proof of Theorem \ref{thm1.5}, choose vectors $\{v_1,...,v_{r_0}\}$ in $V$ and
$\{v_1^*,...,v_{r_0}^*\}$ in
 $V^*$ so $v_i^*\SS(v_0)v_j=\delta_{ij}$. Let $f(v):=\det\{v_i^*\SS(v)v_j\}$. Since
$f(v_0)\ne0$, the holomorphic function $f$ does not vanish identically and thus there exists a timelike
vector $v^-$ in $V$ so
$f(v^-)\ne0$. This proves $r_-\ge r_0$.

To show that $r_0<r_-$, we suppose to the contrary that $r_0=r_+=r_-$ and
argue for a contradiction; this implies $\SS(v)$ has constant rank $r$ on $V-\{0\}$. Let $V_+$ be a maximal space
like subspace of $V$ and let $V_-:=V_+^\perp$ be the complementary maximal timelike subspace. If $v\in V$, then
we may decompose $v=v^++v^-$ for $v^\pm\in V^\pm$. We define a self-adjoint map linear map $\psi$ of $V$ and a positive
definite inner product on $V$ by setting:
\begin{eqnarray}
&&\psi(v):=v^+-v^-,\quad\text{and}\nonumber\\
&&(v,w)_e:=(\psi v,w)=-(v_-,w_-)+(v_+,w_+).\label{eqn5.1}\end{eqnarray}
Define $A(v):=\SS(\psi v)\psi$. We show that $A$ is self-adjoint with respect $(\cdot,\cdot)_e$:
\begin{eqnarray*}
(A(v)v_1,v_2)_e&=&(\SS(\psi v)\psi v_1,v_2)_e=(\psi\SS(\psi v)\psi v_1,v_2)\\
&=&(v_1,\psi\SS(\psi v)\psi
v_2)=(v_1,A(v)v_2)_e.\end{eqnarray*}
Let $S^{p+q-1}:=\{v\in V:(v,v)_e=1\}$ be the associated unit sphere. Because we have that $A(v)v=\SS(\psi v)\psi v=0$, we
have an orthogonal direct sum decomposition
$$T_v(S^{p+q-1})=E_+(v)\oplus E_0(v)\oplus E_-(v),$$
where $E_+(v)$, $E_-(v)$, and $E_0(v)$ are the span of the eigenvectors of $A(v)$ corresponding to positive eigenvalues,
the zero eigenvalue, and negative eigenvalues, respectively. Since
$A(-v)=-A(v)$,
$E_\pm(v)=E_\mp(-v)$. As
$A$ has constant rank, $\dim(E_\pm)$ is constant so these spaces fit together to define smooth bundles on $S^{p+q-1}$.

Since
$\SS\ne0$,
$E_\pm(v)\ne0$. We choose a section $s_+$ to $E_+$ only vanishing at the north pole. Then $s_-(v):=s_+(-v)$ is a section
which only vanishes at the south pole and which satisfies $s_+\perp s_-$. Setting $s:=s_++s_-$ then constructs a nowhere
vanishing section to $T(S^{p+q-1})$ with $s(-v)=s(v)$; this contradicts Lemma \ref{lem2.3} and completes the proof of
Theorem \ref{thm1.6} (1).\end{proof}

The following is a useful technical Lemma.
\begin{lemma}\label{lem5.1} Let $V$ be a vector space of signature $(p,q)$. Let $X$ be a connected topological space on
which $\mathbb{Z}_2$ acts. Let $f$ be a continuous map from $X$ to the space of self-adjoint linear transformations of $V$.
Assume that $r=\rank(f)$ is constant and that $f(-x)=-f(x)$. Then $r$ is even.
\end{lemma}

\begin{proof} We adopt the notation of equation (\ref{eqn5.1}) to define $(\cdot,\cdot)_e$ and $\psi$. We shall let
$\tilde f(x):=f(x)\psi$. Decompose the trivial bundle $X\times V=E_-\oplus E_0\oplus E_+$ as the span of the eigenvectors
of $\tilde f$ corresponding to positive eigenvalues, the zero eigenvalue, and negative eigenvalues, respectively. Since
$f$ has constant rank, these define vector bundles over $X$. Since $E_-(x)=E_+(-x)$ and as $X$ is connected, $\dim E_-=\dim
E_+$. Since $\tilde f$ is self-adjoint with respect to a positive definite inner product, $\tilde f$ vanishes on $E_0$. Thus
$r=\rank(f)=\rank(\tilde f)=\dim E_-\oplus\dim E_+=2\dim E_-$ is even.
\end{proof}

\begin{proof}[Proof of Theorem \ref{thm1.6} (2)] Let $V$ be a vector space of signature $(p,q)$ where $q\equiv2$
mod $4$ and where $p<q-1$. By Theorem \ref{thm1.2}, we may assume $p\ge2$. Let $\RR$ be a Jordan Szab\'o tensor on $V$. We
must show
$\RR=0$. By equation (\ref{eqn1.4}), $\nu(q)=1$. Because $p<q-1=q-\nu(q)$, Theorem
\ref{thm1.4} (2) shows $r_+\le 2$. Thus, by Theorem \ref{thm1.6} (1), $r_0<2$. Since by Lemma \ref{lem5.1} we have $r_0$ is
even, we conclude therefore that $r_0=0$. Choose a basis $e_i$ for $V$ and expand $\SS(x)v_i=\sum_jS_{ij}(x)v_j$. The
functions
$S_{ij}$ are cubic polynomials in $x$ which vanish identically on the real null-cone. Complexification then extends this
relationship to the complex null cone. Since the polynomial $(x,x)$ is irreducible as $\dim(V)\ge3$, we may use the Hilbert
Nullstellungsatz to see there exists a linear function $f(x)$ so that:
$$\SS(x)=(x,x)f(x).$$

Since
$\SS(x)x=0$, we have $f(x)x=0$ for $(x,x)\ne0$ and hence $f(x)x=0$ for all $x$ by continuity. We polarize this relation to
see $f(x)y+f(y)x=0$ for all $x,y\in V$. Dotting this relationship with $y$ then yields
$$0=(f(x)y,y)+(f(y)x,y)=(f(x)y,y)+(x,f(y)y)=(f(x)y,y)$$
for all $x,y\in V$. Polarization then yields $(f(x)y,z)+(f(x)z,y)=0$ for all $x,y,z$ in $V$ and hence, as $f$ is symmetric,
$2(f(x)y,z)=0$ for all $x,y,z$ in $V$. Since the metric on $V$ is non-degenerate, we may conclude that $f(x)$ and hence
$\SS(x)$ vanishes identically. Theorem \ref{thm1.6} (2) now follows from Lemma \ref{lem3.3}. \end{proof}


\begin{thebibliography}{AAA}



\bibitem{refAdams} J. Adams, {\it Vector fields on spheres}, Annals of Math., {\bf 75},
(1962), 603--632.

\bibitem{refFG02} B. Fiedler and P. Gilkey, {\it Nilpotent Szab\'o, Osserman and Ivanov-Petrova pseudo-Riemannian
manifolds}, (preprint) http://arXiv.org/abs/math.DG/0211080.

\bibitem{refGRKVL} E. Garci\'a-Ri\'o, D. Kupeli, and R. V\'azquez-Lorenzo, {\bf Osserman Manifolds in Semi-Riemannian
Geometry}, Lecture notes in Mathematics, 1777, Springer Verlag (2002), ISBN 3-540-43144-6.



\bibitem{refGil} P. Gilkey, {\bf Geometric Properties of Natural Operators Defined by the Riemann
Curvature Tensor}, World Scientific Press (2001), ISBN 981-02-04752-4.

\bibitem{refGiIv} P. Gilkey and R. Ivanova, {\it Spacelike Jordan Osserman algebraic curvature tensors in
signature $(p,q)$ for
$p<q$}, (preprint: http://arXiv.org/abs/math.DG/0205072).

\bibitem{refGIZ} P. Gilkey, R. Ivanova, and T. Zhang, {\it Szab\'o Osserman IP Pseudo-Riemannian manifolds},
(preprint: http://arXiv.org/abs/math.DG/0205085).

\bibitem{refGilStav} P. Gilkey and I. Stavrov, {\it Curvature tensors
   whose Jacobi or Szab\'o operator is nilpotent
    on null vectors}, Bull. London Math. Soc. (to appear);\newline http://arXiv.org/abs/math.DG/0205074.

\bibitem{refGiZa} P. Gilkey and T. Zhang, {\it Algebraic curvature tensors for indefinite metrics whose
skew-symmetric curvature operator has constant Jordan normal 
form}, Houston Math J. {\bf 28} (2002), 311--328.

\bibitem{refIP}  S. Ivanov and I. Petrova,  {\it Riemannian manifold in which the
skew-symmetric curvature operator has pointwise constant eigenvalues},  Geom. Dedicata,
{\bf 70}, (1998), 269--282.

\bibitem{refRIGS}  R. Ivanova and G. Stanilov, {\it A skew-symmetric curvature
operator in Riemannian geometry}, in {\bf  Sympos. Gaussiana, Conf A}, ed. Behara, Fritsch,
and Lintz (1995), 391--395.


\bibitem{refOss} R. Osserman, Curvature in the eighties, {\it Amer. Math.
    Monthly}, {\bf 97}, (1990), 731--756.

\bibitem{refSpanier} E. Spanier, {\bf Algebraic Topology}, McGraw-Hill (1966); corrected reprint Springer-Verlag (1981) ISBN:
0-387-90646-0.

\bibitem{refSt} I. Stavrov, Ph.D. Thesis, University of Oregon (2003).

\bibitem{refSzabo}
Z. I. Szab\'o, {\it A short topological proof for the symmetry of $2$ point homogeneous
spaces}, Invent. Math., {\bf 106}, (1991), 61--64.


\bibitem{refWolf} J. Wolf, {\bf Spaces of constant curvature,} 5th ed., Publish or Perish Inc. (1984), ISBN
0-914098-07-1

\end{thebibliography}
\end{document}